\documentclass
[12pt,a4paper]{article}
\usepackage[cp1251]{inputenc}
\usepackage{amssymb, amsfonts, amsmath, latexsym}
\usepackage[russian]{babel}
\begin{document}

\Large
\centerline{\bf Box Resolvability}\vspace{6 mm}

\normalsize\centerline{\bf  Igor Protasov}\vspace{6 mm}

{\bf Abstract.} We say that a topological group $G$  is partially box $\kappa$-resolvable if there exist a dense  subset $B$ of $G$  and a subset $A $ of $G$, $|A|=\kappa$ such that the subsets $\{ aB: a\in A\}$ are pairwise disjoint. If $G=AB$ then $G$ is called box $\kappa$-resolvable. We prove two theorems. If a topological  group $G$ contains an injective convergent sequence then $G$ is box $\omega$-resolvable. Every infinite totally bounded topological group $G$ is partially box $n$-resolvable for each natural number $n$, and $G$ is box $\kappa$-resolvable for each infinite cardinal $\kappa, \kappa<|G|$.


\vspace{6 mm}

\centerline{\bf 1.	Introduction}

\vspace{6 mm}

For a cardinal $\kappa$, a topological space $X$ is called {\it $\kappa$-resolvable} if $X$ can be partitioned into $\kappa$ dense subsets [1].  In the case $\kappa=2$, these spaces were defined by Hewitt [6] as {\it resolvable spaces}. If $X$ is not $\kappa$-resolvable then $X$ is called {\it $\kappa$-irresolvable.}

In topological groups, the intensive study of resolvability was initiated by the following remarkable theorem of Comfort and van Mill [3]:
every countable  non-discrete Abelian topological group $G$ with finite subgroup $B(G)$ of elements  of order 2 is 2-resolvable.  In fact [19], every infinite Abelian  group $G$ with finite $B(G)$ can be partitioned into $ \omega$ subsets dense in every non-discrete group topology on $G$. On the other hand, under MA, the countable Boolean group $G$, $G=B(G)$ admits maximal (hence, 2-irresolvable)  group topology [8].  Every non-discrete $\omega$-irresolvable topological group $G$ contains an open countable Boolean subgroup provided  that $G$ is Abelian [11] or countable [18], but the existence of non-discrete $\omega$-irresolvable group topology on the countable Boolean group implies that there is a $P$-point in $\omega^{\ast}$ [11].  Thus, in some models of ZFC (see [14]), every    non-discrete Abelian or countable topological group is $\omega$-resolvable. We mention also $\kappa$-resolvability of every infinite totally bounded topological group $G$ of cardinality $\kappa$ [9]. For systematic exposition of resolvability in topological and left topological group see [4, Chapter 13].

This note is to introduce more delicate kind of resolvability, the box resolvability.

Given a group $G$  and a cardinal $\kappa$, we say that a subset $B$ of $G$ is a {\it partial box} of index $\kappa$, if there exists a subset $A$ of $G$,  $|A|= \kappa$ such that the subsets $\{aB: a\in A\}$ are pairwise disjoint. In addition, if $G=AB$ then $B$ is called a {\it box} of index $\kappa$. Example: a subgroup $H$ of $G$ is a box of index $|G: H|$, and any set $R$  of representatives of right cosets of $G$ by $H$ is a box of index $|H|$.

We use also the factorization terminology [16]. For subset $A,B$  of $G$, the product $AB$ is called a {\it partial factorization} if $aB \bigcap a^{\prime}B=\emptyset$  for any distinct $a,a^{\prime}\in A$. If $G=AB$ then the product $AB$ is called a {\it factorization} of $G$.   Thus, a box $B$ of index  $\kappa$ is a right factor of some factorization $G=AB$  such that $|A|=\kappa$.

We say that a topological group $G$ is {\it (partially) box} $\kappa$-resolvable if there exists a (partial) box $B$ of index $\kappa$ dense in $G$. Clearly, every partially box $\kappa$-resolvable topological group is  $\kappa$-resolvable, but a $\kappa$-resolvable group needs not to be box $\kappa$-resolvable (see Examples 1 and 2). However, I do not know, whether every 2-resolvable group is partially box 2-resolvable.

On exposition: in section 2, we prove two theorems announced in Abstract and discuss some prospects of box resolvability in section 3.
\vspace{6 mm}

\centerline{\bf 2.	Results}
\vspace{6 mm}

We begin with two examples demonstrating purely algebraic obstacles to finite box resolvability.

\vspace{3 mm}
{\bf Example 1. } We assume that the group $\mathbb{Z}$ of integer numbers is factorized $\mathbb{Z}=A+B$ so that $A$ is finite, $|A|>1$. By the Haj$\acute{o}$s theorem [5], $B$ is periodic: $B=m+B$ for some $m\neq 0$. Then $m\mathbb{Z} +B=B$  and $m\mathbb{Z}+ b\subseteq B$ for $b\in B$.

Now we endow $\mathbb{Z}$ with the topology $\tau$ of finite indices (having $\{ n\mathbb{Z}: z\in \mathbb{N} \}$ as the base at 0). Since $m \mathbb{Z}$ is open in $\tau$ and $|A|>1$, we conclude that $B$ is not dense,  so ($\mathbb{Z},\tau $) is box $n$-irresolvable for each $n>1$. $\Box$
\vspace{3 mm}

{\bf Example 2. } Every torsion group $G$ without elements of order 2 has no boxes of index 2. We assume the contrary: $G=B\bigcup gB$, $B\bigcap gB=\emptyset$  and $e\in B$, $e$ is the identity of $G$.  Then $g^{2}B=B$  and $B$ contains the subgroup $<g^{2}>$ generated by $g$. Since $g$ is an element of odd order, we have $g\in <g^{2}>$ and $g\in B\bigcap gB$.$\Box$
\vspace{3 mm}

Let $G$ be a countable group. Applying [13, Theorem 2], we can find a factorization $G=AB$  such that $|A|=|B|=\omega.$ Hence, if we endow $G$ with a group topology $\tau$, there are no algebraic obstacles to box $\omega$-resolvability of $(G,\tau)$. \vspace{3 mm}

In what follows, we use two elementary observations. Let $G$ be a topological group, $H$ be a subgroup $G$ and $R$ be a system of representatives of right cosets of  $G$ by $H$.  Let $AB$ be a factorization of $H$. Then we have

$(1)$	If $B$ is dense in $A$ then $A(BR)$ is a factorization of $G$ with dense $BR$;

$(2)$	If $R$ is dense in $G$ then $A(BR)$ is a factorization of $G$ with dense $BR$.
\vspace{6 mm}

{\bf Example 3. } Let $G$ be a non-discrete metrizable group and let $A$ be a subgroup of $G$. If $A$ is either finite or countable discrete then there is a factorization $ AB$ of $G$ such that $B$ is dense in $G$.

In view of (1), we may suppose that $G$ is countable. Let $\{U_{n}: n\in\omega \}$ be a base of topology of $G$. For each $n\in\omega$, we choose $x_{n}\in U_{n}$ so that  $Ax_{n}\bigcap Ax_{m} =\emptyset$ if $n\neq m$. Then we complement the set $\{x_{n}: n\in\omega\}$ to  some full system $B$ of representatives of right cosets  of $G$ by $A$.   $\square$\vspace{3 mm}

If a topological group $G$ contains an injective convergent sequence then $G$ is  $\omega$-resolvable (see [2, Lemma  5.4]). If an injective sequence $(a_{n})_{n\in\omega }$ converges to the identity  $e$ in some group topology on a group $G$   then, by [13, Theorem 1], the set $\{e, a_{n}, a_{n} ^{-1} : n\in\omega\}$ is a left factor of some factorization of $G$.
\vspace{6 mm}

{\bf Theorem 1} {\it If a topological group $G$ contains an injective convergent sequence $(a_{n})_{n\in\omega }$  then $G$ is box $\omega$-resolvable.
\vspace{3 mm}

Proof.} We suppose that $(a_{n})_{n\in\omega }$  converges to the identity $e$ of $G$ and denote
$$A=\{ e, a_{n}, a_{n} ^{-1} : n\in\omega \},  \\   A_{n}=\{ e, a_{m}, a_{m} ^{-1} : m\leq n\},  C_{n}= A\setminus A _{n} .$$

Replacing $G$ to the subgroup of $G$ generated by $A$, in view of (1), we may suppose that $G$ is countable, $G=\{ g_{n}: n\in\omega \}$,  $g_{0}=e.$  Our goal is to find a factorization $G=AB$ such that $B$ is dense in $G$.  We shall construct a family $\{B_{n}: n\in\omega\}$, $B_{n}\subset B_{n+1}$  of finite subsets of $G$ such that, for each $n\in\omega $,

$(3)$	$AB_{n}$ is a partial factorization;

$(4)$	$\{g_{0}, …, g_{n}\}\subset AB_{n};$

$(5)$	for every $g \in A_{n-1}  B_{n-1}$, there exists $b\in B_{n}$ such that $g\in C_{n} ^{2} b.$

After $\omega$ steps, we put $B=\bigcup _{n\in\omega} B_{n} .$  By (3)  and (4), $AB$ is a factorization  of $G.$
By (4)  and (5), $B$ is dense in $G$.

We put $B_{0}=\{ e\}$ and suppose that we have chosen $B_{0}, …, B_{n}$ satisfying (3) , (4) and (5).

To make the inductive step from $n$ to $n+1$, we use the following observation.

$(6)$  If $F$ is a finite subset of $G$ and $g\notin AF$ then there is $k\in\omega $ such that $AC _{k}g \bigcap AF =\emptyset$.
\vspace{3 mm}

Indeed, if  $AC _{k}g \bigcap AF =\emptyset$ for each $k\in\omega$ then there are an injective sequence $(s_{n})_{n\in\omega}$  in $A$ and $a\in A$ such that $as_{n} g\in  AF$  for each $n$, so $g\in a^{-1} F$  and $g\in AF$.

We choose the first element $g\in\{ g_{n}: n\in\omega\}\setminus AB_{n}$ and use (6) with $F=B_{n}$ to find $k\in\omega$ such that $AC _{k} g\bigcap AB_{n}=\emptyset .$   We take $c\in C _{k} $ and notice that $Acg\bigcap AB_{n} =\emptyset$ and $g\in Acg.$

We enumerate $x_{0}, …, x_{p}$ the elements of the set $A_{n} B_{n}\setminus B_{n}$ and take $s\in C _{n+1}$ such that, for each $i\in\{ 0,…, p\}$,
$$sx_{i}\notin  A_{n}B_{n} \bigcup Acg.$$
Then we use (6)  to choose $c_{0},…, c_{p} \in C _{n+1}$ such that,  for each $i\in\{0,..,p\}$, $Ac_{i} (sx_{i})\bigcap $
$\{sx_{i+1}, ..., sx_{p}\}=\emptyset $ and $$Ac_{j}(s x_{i})\bigcap (AB_{n}\bigcup Acg\bigcup Ac_{0}(sx_{0})
\bigcup ...\bigcup Ac_{i-1}(sx_{i-1}))=\emptyset .$$

After that, we put $$B _{n+1} = B_{n} \bigcup\{cg, c_{0}s x_{0}, ..., c_{p} s x_{p}\} $$ and note that (3), (4),  (5) hold for $n+1$ in place of $n$.
$\square$

\vspace{6 mm}

{\bf Theorem 2.} {\it Let $G$ be an infinite totally bounded topological group of cardinality $\gamma$, $H$ be a subgroup of $G$ such that $|G: H|=\gamma$, $F$ be a finite subset of $G$. Then the following statements hold

$(i)$ there is a partial factorization $FB$ such that $B$ is dense in $G$;

$(ii)$ there is a factorization $G=HR$ such that $R$ is dense in $G$;

In particular, $G$ is a box $n$-resolvable for each $n\in \mathbb{N}$, and $G$ is box $\kappa$-resolvable for each infinite cardinal  $\kappa$, $\kappa<\gamma$.\vspace{3 mm}

Proof . } $(i)$ We denote $\mathfrak{F} _{F}= \{ K\subset G: |K|<\omega ,  K^{-1} K \bigcap F^{-1} F=e\}$ and enumerate $\mathfrak{F} _{F}= \{ K_{\alpha}: \alpha<\gamma\}$. We choose inductively a $\gamma$-sequence $(x_{\alpha})_{\alpha< \gamma} $ in $G$ such that

\vspace{3 mm}
(7) $FK_{\alpha}^{-1} x _{\alpha}\bigcap F K_{\beta}^{-1}x _{\beta}=\emptyset$ for all $\alpha,\beta,$ $\alpha<\beta<\gamma,$ and denote $B=\bigcup_{\alpha<\gamma} K_{\alpha}^{-1}x _{\alpha}$. We take distinct $g, h\in F$. Since $K_{\alpha}\in \mathfrak{F}_{H}$, we have $g K_{\alpha}^{-1}x _{\alpha}\bigcap h K_{\alpha}^{-1}x _{\alpha}=\emptyset $. If $\alpha<\beta$  then  , by (7), $g K_{\alpha}^{-1}x _{\alpha}\bigcap h K_{\beta}^{-1}x _{\beta}=\emptyset $. Hence, $gB\bigcap hB = \emptyset$ and the product $FB$ is a partial factorization.
\vspace{3 mm}

To prove that $B$ is dense, we use

(8) for any open subsets $U_{1}, …, U_{n}$ of $G$, there exist $y_{1}\in U_{1}, …, y_{n}\in  U_{n}$  such that $\{y_{1},…,y_{n}\}\in  \mathfrak{F}_{F}$,  that can be easily proved by induction on $n$.

Now let $U$ be a neighborhood of $e$ and $g\in G$.  We show that $Ug\bigcap B\neq\emptyset.$  We take a  neighborhood $V$ of $e$ such that $V^{-1}V\subseteq U$. Since $G$ is totally bounded, there are $z_{1},…, z _{n }\in G$  such that $G= z_{1} V\bigcup ... \bigcup z_{n}V$. We use (8) to find $y_{1}\in z_{1} V,…, y_{n} \in z_{n} V$ such that $\{y_{1},…, y_{n}\}\in \mathfrak{F}_{F}$. Then $z_{1}\in y_{1} V^{-1}, …, z_{n} \in  y_{n} V^{-1}$ so  $\{y_{1}, …, y_{n}\}  U=G$. We chose $\alpha<\gamma $ such that $K_{\alpha} = \{y_{1}, …, y_{n}\} $. Since $K_{\alpha}U_{g}=G$,  we have $x_{\alpha}\in  K_{\alpha}Ug$, $K_{\alpha}^{-1} x_{\alpha} \bigcap Ug\neq\emptyset$ and $B\bigcap Ug\neq\emptyset$.

\vspace{3 mm}

$(ii)$ For any open    subset $U$ of $G$,  we choose a finite subset $F_{U}$ such that $G=F^{-1} _{U} U$ and $Hx\bigcap Hy=\emptyset$  for all distinct $x,y\in  F_{U}$. We enumerate without repetitions the set $\{F_{U}: U$ is open $\}$ as $\{K_{\alpha} : \alpha<\gamma\}$. Since $|G: H|=\gamma$,  we can choose inductively a $\gamma$-sequence $(x_{\alpha})_{\alpha<\gamma}$ in $G$  such that $HK_{\alpha}x_{\alpha}\bigcap H K_{\beta} x_{\beta}=\emptyset$ for all $\alpha<\beta<\gamma.$ We denote $S=\bigcup_{\alpha<\gamma}K_{\alpha} x_{\alpha}$ and show that $S$ is dense in $G$. Given any open subset $U$ in $G$, we choose $\alpha<\gamma$ such that $G=K^{-1}_{\alpha} U$. Then $x_{\alpha}\in K^{-1}_{\alpha} U$ so $K_{\alpha}  x_{\alpha}\bigcap U\neq\emptyset$ and $S\bigcap U \neq\emptyset$.

To conclude the proof, we complement $S$ to some full system $R$ of representatives of right cosets of $G$ by $H$.
$\Box$

\vspace{6 mm}

\centerline{\bf 3.	Comments}
\vspace{6 mm}

1.	In connection with Theorem  2, we should ask
\vspace{3 mm}

{\bf Question 1.} {\it Is every infinite totally bounded topological group of cardinality $\kappa$ box  $\kappa$-resolvable?}\vspace{3 mm}

For  $\kappa=\omega$ , to answer this question positively, it suffices to generalize Theorem 1 and prove that a topological  group $G$ is box $\omega$-resolvable provided that $G$ contains a countable thin subset $X$ such that $e$ is the unique limit point of $X$. A subset $X$ of $G$ is called {\it thin} if $|gX \bigcap X|<\omega$  for every $g\in G\setminus  \{e\}$. By [12, Theorem 2], every infinite totally bounded topological group $G$ have such a subset $X$.

In the case $\alpha=\omega$, I believe in the  positive answer to Question 1 but then
\vspace{3 mm}

{\bf Question 2.} {\it In ZFC,  does there exist an infinite non-discrete box $\omega$-irresolvable topological group? }
\vspace{3 mm}

2.	Given a family $\mathcal{I}$ of subsets of a topological group $G$, we say that $G$ is {\it $\mathcal{I}$-box  $\kappa$-resolvable } if there exist a dense subset $B$  of $G$ and a subset $A$ of cardinality $\kappa$  such that $G=AB$  and $aB \bigcap a' B\in \mathcal{I}$ for all distinct $a,a' \in A$.

{\it Every countable totally bounded topological group  $G$ is $\mathcal{I}$-box $\omega$-resolvable with respect to the family $\mathcal{I}$ of all finite subsets of $G$.}

By [12, Theorem 3], $G$ has a thin dense subset $B$. Then $gB\bigcap  g'B\in I$  for all distinct $g, g'\in G$ and $G=GB$.
\vspace{3 mm}

{\bf Question 3.} Let $\tau$ be the topology of finite indices (see Example 1) on $\mathbb{Z}$. Is $(\mathbb{Z},\tau)$ $\mathcal{I}$-box 2-resolvable with respect to the family $\mathcal{I}$ of all nowhere dense subsets of $\mathbb{Z}$?

\vspace{3 mm}

3.	The notion of the box resolvability is natural in much more general context of $G$-spaces. Let $X$ be a topological space and let $G$ be a discrete group. We suppose that $X$ is endowed with transitive action $G\times X\rightarrow  X : (g,x)\longmapsto  gx$  such that, for each $g\in  G$, the mapping $x\longmapsto gx$ is continuous.

We say that $X$  is {\it box $\kappa$-resolvable } if  there exist a dense subset $B$ of $X$  and a subset $A$  of $G$, $ |A|=\kappa $  such that $X=AB$  and the subsets $\{aB:  a\in A\}$ are pairwise disjoint.

For example, take the group $\mathbb{Q}$ of rational number with the natural topology and let $G$ be a group of all homomorphisms of   $X$  such that, for   each $g\in G$, there is $a\in \mathbb{Q}$, $a>0$ such that $gx=x$  for every $x\in \mathbb{Q} \setminus [-a, a]$. Then $\mathbb{Q}$  is box $\kappa$-resolvable only for $\kappa=1$ and Theorem 1 does not hold for some $G$-spaces. On the other side if $G$ is the group of all homeomorphisms of $\mathbb{Q}$  then, by Theorem 1, $\mathbb{Q}$ is box  $\omega$-resolvable because $G$ contains the subgroup of translations of $\mathbb{Q}$.
\vspace{3 mm}

4.	A topology on a group $G$ is called {\it left invariant} if all shifts $x\mapsto  gx,$ $g\in G$ is continuous, and a group $G$ endowed with a left invariant topology is called  {\it left topological}. Clearly, every left topological group has the natural structure of $G$-space.

In ZFC, every infinite group  $G$ admits maximal (hence, irresolvable) regular left invariant topology [10].

Every non-discrete left topological group of second category is  $\omega$-resolvable [4, Theorem 13.1.12], but in some model of ZFC there is  an irresolvable homogeneous space of second category [15].
\vspace{3 mm}

{\bf Question 4.} {\it Is every box $\omega$-irresolvable left topological group meager? }
\vspace{3 mm}

5.	We say that a left topological group is {\it locally box $\kappa$-resolvable} if there exists  a subset  $B$ of $G$ such that, for each neighborhood $U$ of the identity $e$, we can choose $A\subseteq U$ such that $|A| =\kappa $, $e \in A$,  $AB$ is a partial factorization and the closure of each subset $aB$, $a\in A$  is a neighborhood  of $e$.

If $G$ is locally box 2-resolvable then some  neighborhood  of $e$ can be partitioned into two dense subsets so $G$ is 2-resolvable.

To see that the converse statement does not hold, we can use the semigroup structure in the Stone-$\check{C}$ech compactification of a discrete group (see [7]). Given an infinite group $G$, we choose  two idempotents $p$ and $q$ from $\beta G\setminus G$ such that $pq=q$, $qp=p$.  We take the family $\{P\bigcup Q\bigcup \{e\} : P\in p, Q\in q$  as the base at $e$ for some left invariant topology $\tau$.
Then $(G,\tau)$ is 2-resolvable but locally box  2-irresolvable. Moreover, under MA, on a countable Boolean group $G$, there $p,q$ such that corresponding $\tau$  is a group topology [17].
\vspace{3 mm}

{\bf Question 5.} {\it In ZFC, does there exists a locally box 2-irresolvable topological group? }
\vspace{3 mm}

\centerline{\bf References}\vspace{3 mm}

[1]  J. Ceder, On maximally resolvable spaces, Fund.Math. 55 (1964) 87-93.

[2] W. Comfort, S. Garcia-Ferreira, Resolvability: a selective survey and some new results, Topology Appl. 74 (1996) 149-167.

[3] W. Comfort, J. van Mill, Group with only resolvable group topologies, Proc. Amer. Math. Soc. 120(1993) 687-696.

[4] M. Filali, I. Protasov, Ultrafilters and Topologies on  Groups, Math. Stud. Monogr. Ser., vol 13, VNTL, Lviv, 2010.

[5] G. Haj$\acute{o}$s, Sur la factorization des groups abelianes, $\check{C}$asopis Pest. Mat. Fys. 74(1949) 157-162.

 [6] E. Hewitt, A problem of set-theoretic topology, Duke Math. J. 10 (1943) 309-333.

[7] N. Hindman, D. Strauss, Algebra in the Stone-$\check{C}$ech Compactification, 2nd edition, de Gruyter, 2012.

[8] V. Malykhin, Extremally disconnected and similar groups, Soviet Math. Dokl. 16 (1975) 21-25; translation from Dokl. Akad. Nauk SSSR 220 (1975) 27-30.

[9]  V. Malykhin, I. Protasov, Maximal resolvability of bounded groups, Topology Appl. 73 (1996) 227-232.

[10] I. Protasov, Maximal topologies on groups, Siberian Math. J.  39 (1998) 1184-1194; translation from Sibirsk. Mat. Zh. 39 (1998) 1368-1381.

[11] I. Protasov, Irresolvable topologies on groups, Ukr. Math. J. 50 (1998) 1879-1887; translation from Ukr. Math. Zh. 50 (1998) 1646-1655.

[12] I. Protasov, Thin subsets of topological group, Topology Appl. 160(2013) 1083-1087.

[13] I. Protasov, S. Slobodianik, Kaleidoscopical configurations in groups, Math. Stud. 36 (2011) 115-118.

[14] S. Shelah, Proper Forcing, Lecture Notes Math, vol. 940, Springer-Verlag, 1982.

[15] S. Shelah, Baer irresolvable spaces and lifting for the layered ideals, Topology Appl. 33 (1989) 217-221.

[16] S. Szab$\acute{o}$, A Sanders, Factoring Groups into Subsets, CRC Press, 2009.

[17] E. Zelenyuk, Topological groups with finite semigroups of  ultrafilters, Mat. Stud. 6 (1996)  41-52.

[18] E. Zelenyuk, On partition of groups into dense subsets, Topology Appl. 126 (2000) 327-339.

[19] E. Zelenyuk, Partitions and sums with inverses in Abelian groups, J. Combin. Theory, Ser. A 115 (2008) 331-339.
\vspace{6 mm}

Department of Cybernetics, Kyiv University,
Prospect Glushkova 2, corp. 6,
03680 Kyiv, Ukraine
e-mail: I.V. Protasov@gmail.com

\end{document}